\documentstyle[12pt,epsf]{article}
\begin{document}
\title{Singularities
\footnote{A chapter for the book ``History of Topology'', ed I. M. James} }
\author{Alan H. Durfee
\footnote{Author address:  Department of Mathematics, Statistics and
Computer Science, Mount Holyoke College, South Hadley, MA 01075 USA.
email: adurfee@mtholyoke.edu.  web site:  http://www.mtholyoke.edu/\~{}adurfee}}

%%%%%%%%%%%%%%%%%%%%%%%%%%%%%%%

%postscript macros from Fishwick article
\newcommand{\postscript}[2]
{\setlength{\epsfxsize}{#2\hsize}
\centerline{\epsfbox{#1}}}

%%%%%%%%%%%%%%%%%%%%%%%%%%%%%%

\def\mod4{\ (\mbox{mod \ }4)}
\def\real{{\bf R}}
\def\complex{{\bf C}}
\def\projective{{\bf P}}
\def\sphere{{\bf S}}
\def\ball{{\bf D}}
\def\origin{{\bf 0}}
\def\x{{\bf x}}
\def\p{{\bf p}}
\def\boundary{\partial}
\def\infinity{\infty}
\def\del{\partial}
\def\euler{\chi}
\def\overbar{\overline}
\def\intersect{\cdot}
\def\congruent{\equiv}

%%%%%%%%%%%%%%%%%%%%%%%%%%%
\def\umlaut{\"}
\def\grave{\`}
\def\aigu{\'}

%%%%%%%%%%%%%%%%%%%%%%%%%%%%%%%%%%%

%numbering commands

\newcommand{\mysection}[1]{\setcounter{mysubsection}{0}\section{#1}}
\newcounter{mycounter}[section]
\renewcommand{\themycounter}{\arabic{section}.\arabic{mycounter}}

%The environments for theorem, proposition, lemma, corollary
%get numbers and italics

\newenvironment{theorem}%
{\medskip 
    \refstepcounter{mycounter}
    {\bf \noindent Theorem.}\begin{em}}%
{\end{em} \medskip }

\newenvironment{proposition}%
{\medskip 
    \refstepcounter{mycounter}
    {\bf \noindent Proposition. \ } \begin{em} }%
{\end{em} \medskip }

\newenvironment{lemma}%
{\medskip 
    \refstepcounter{mycounter}
    {\bf \noindent Lemma \themycounter. \ } \begin{em} }%
{\end{em} \medskip }

\newenvironment{corollary}%
{\medskip 
    \refstepcounter{mycounter}
    {\bf \noindent Corollary. \ } \begin{em} }%
{\end{em} \medskip }

\newenvironment{formula}%
{\medskip 
    \refstepcounter{mycounter}
    {\bf \noindent Formula \themycounter. \ } \begin{em} }%
{\end{em} \medskip }

%The environments for remark, definition, example, conjecture and problem get numbers
%but not italics. 

\newenvironment{remark}%
{\medskip 
    \refstepcounter{mycounter}
    {\bf \noindent Remark \themycounter. \ }}%
{\medskip }

\newenvironment{conjecture}%
{\medskip 
    \refstepcounter{mycounter}
    {\bf \noindent Conjecture. \ }}%
{\medskip }

\newenvironment{definition}%
{\medskip 
    \refstepcounter{mycounter}
    {\bf \noindent Definition \themycounter. \ }}%
{\medskip }

\newenvironment{example}%
{\medskip 
    \refstepcounter{mycounter}
    {\bf \noindent Example \themycounter. \ }}%
{\medskip }

\newenvironment{problem}%
{\medskip 
    \refstepcounter{mycounter}
    {\bf \noindent Problem \themycounter. \ }}%
{\medskip }

%The proof environment gets neither italics nor numbers, 
%but gets an ``end of proof sign''.

\newenvironment{xproof}%
{\medskip 
  \noindent
    {\bf Proof. \ }}%
{\endofproof \medskip }

%the following command is still in the paper

\newcommand{\resetassertioncounter}{}

%%%%%%%%%%%%%%%%%%%%%%%%%%%%%%%%%%%

\maketitle

\section{Introduction}

This article recounts the rather wonderful interaction of topology and
singularity theory which began to flower in the 1960's with the work of
Hirzebruch, Brieskorn, Milnor and others.
This interaction can be traced back to the work
of Klein, Lefschetz and Picard, and also to the work
of knot theorists at the beginning of this century.
It continues to the present day, flourishing and
expanding in many directions.
However, this is not a survey article, but a history; the events of our
time are harder to see in perspective, harder to marshal
into coherent order, and their very multitude makes it impossible to
recount them all.
Hence this interaction is followed forward in only a few
directions.
\footnote{That I have attempted to do this at all is due to the prodding
of my conscience and a list suggested by W. Neumann of some recent
areas where topology has had an effect on singularity theory.  He added,
though, that ``the task becomes immense...other people would probably
come up with almost disjoint lists.''  
%I have been adding names and
%topics until the day this article had to sent off for publication, so 
The
randomness of my efforts here should be readily apparent, and my apologies
to those whose work is not mentioned.}

The reader may get a sense of the current state of affairs
in singularity theory by browsing in the conference
proceedings \cite{pspm-40,trieste}.
The focus of this article is singularities of complex algebraic
varieties.
Real varieties are omitted.  
Also omitted from this account is the area of critical
points of differentiable functions, work initiated by Thom, Mather,
Arnold and others; a survey of this subject can be found in the books
\cite{agv-1, agv-2, encyc}.

%In fact, singularity theory is the
%intersection of all branches of mathematics (as asserted by a recent
%colloquium speaker) and hence a large subject.
%Thus some material had to be omitted from a article of this type.

When two areas interact, ideas flow in both directions.
Ideas from topology have
entered singularity theory, where algebraic problems have been understood
as topological problems and solved by topological methods.  
(In fact, often the
crudest invariants of an algebraic situation are topological.)
Conversely, ideas of singularity theory have traveled in the reverse
direction into topology.  Algebraic geometry supplies many interesting
examples both easily and not so easily understood, and these provide a
convenient testing ground for topological theories.

\section{Knots and singularities of plane curves}

In the 1920's and 30's there was much activity in knot theory as the new tools of algebraic topology were being applied;
the fundamental group of the knot complement was introduced, as were the
Alexander polynomial, branched cyclic covers, the Seifert surface,
braids, the quadratic form of a knot, linking invariants, and so forth.
Many clearly-written wonderful papers were produced on these subjects.

At the same time in algebraic geometry there was interest in
understanding complex algebraic surfaces, in particular by exhibiting
them as branched covers of the plane.
This method is the analogous to the method in one dimension lower of projecting a curve to a
line.
The discriminant locus in the latter case is a set of points and it is easy to
understand the branching.  For surfaces the branching is more
complicated since the discriminant locus is a curve.

A method of examining the branching problem for surfaces was proposed by Wirtinger in Vienna,
who gave some seminars on this subject beginning in 1905.
He divided branch points into two types:
At a smooth point of the discriminant curve, the branching group
(``Verzweigungsgruppe'') of the
surface is cyclic, like that of a curve.  
These points were called ``branch points
of type I''.
Singular points of the discriminant curve were
called ``branch points of type II''.
He also worked out a simple example.

The classification and the example were recorded by his student Brauner
in the
beginning of his paper ``On the geometry of functions of two complex
variables'' \cite{brauner}.
Wirtinger's
example is the smooth surface in $\complex^3$ given by the equation
$$z^3 -3zx+2y = 0$$
When this is projected to the $(x,y)$-plane, the discriminant curve is
$$x^3-y^2 = 0$$
There is one point in the surface over the origin in the $(x,y)$-plane,
two points over the remaining points of the curve $x^3 - y^2 = 0$, and
three points over the rest of the plane.

\begin{figure}
\postscript{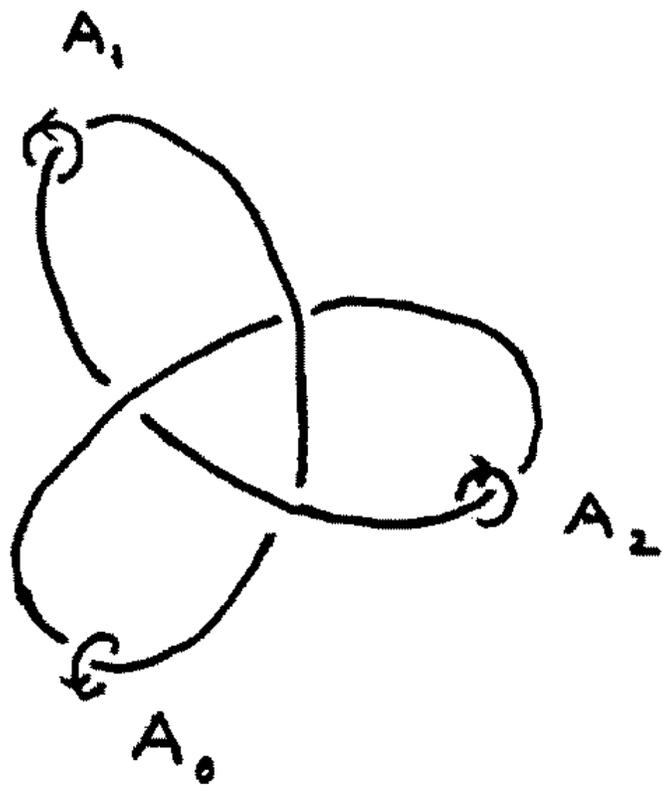}{0.7}
\caption{The discriminant intersected with the sphere}
\label{figure1}
\end{figure}

To understand the type II branching of the surface near the origin, a
three-sphere $\sphere^3_r$ of radius $r$ about the origin in the plane was mapped to real
three-space by stereographic projection.  The image of the intersection
of this three-sphere with the discriminant curve was then exhibited as
a trefoil knot $\Gamma$ (Figure \ref{figure1}).
It sufficed to understand the branching of the surface over $\Gamma \subset
\sphere^3_r$.
Let $A_i$ for $i = 1,2,3$ be the branching substitution produced by
travelling around the loop labelled $A_i$ in the figure.
The $A_i$ must satisfy the (now well-known) Wirtinger relation
$$A_0^{-1}A_1A_0A_2^{-1} = 1$$ 
at the left-hand crossing point of the knot projection in the figure.
The only possibility for the permutation of the sheets of the covering is thus $A_0 = (12)$, $A_1 = (23)$ and $A_2 =
(13)$.
Hence the branching group in the neighborhood of $(0,0)$ is not cyclic (as it
is for plane curves), but the symmetric group on three elements.

Brauner concluded ``Wir haben aus obigem erkannt, dass es die
topologischen Verh\"altnisse der Kurve $\Gamma$ sind, welche dieses
merkw\"urdiger Verhalten der Funktionen in der Umgebung der
Verzweigungsstellen II. Art bedingen.''  [We thus have learned that the
topological form of the curve $\Gamma$ determines this remarkable
behavior of the function in the neighborhood of a type II branch
point.]

There are thus two problems, he said.  The first is to determine the topology of
the (discriminant) curve in the neighborhood of a singular point,
i.e. the knot $\Gamma$. The second is to determine the group given
by the Wirtinger relations (in modern terminology, the fundamental group
of the complement of the knot $\Gamma$).  These two problems were solved in his paper.
He remarked that there are three more problems.
The first is to determine the branching group of a
function locally in the neighborhood of a point. (This group is of
course a quotient of fundamental group of the complement of the knot $\Gamma$.)  Next, one
should determine the global branching group of a function. 
Finally, given a group, is there a function which has this group as
branching group?
 These problems, he said, will form the subject of two further papers
\footnote{I do not know if these papers appeared.}.

He then continued with a systematic study of the links of curve singularities
and their fundamental groups.
He first looked at the curve
$$ax^n + by^m = 0$$
with the g.c.d $(n,m) = 1$,
parameterizing it by setting $x = \alpha t^m$ and
$y = \beta t^n$, where $\alpha$ and $\beta$ were suitably chosen
constants.
He wrote the complex number $t$ as $\rho e^{i\phi}$ with $\rho$ and
$\phi$ real and worked out parametric equations for the intersection of
the curve with the sphere.
Taking its image under the equations for stereographic projection, he
observed that the image curve lay on a torus, winding $n$ times in the
direction of the meridian and $m$ times in the direction of the equator,
and hence was a torus knot.

He then went on to look at two such curves as above and described their
linking. 
He then examined the curve parameterized by 
$x = t^m(a_m + ta_{m+1} + \dots)$ and $y = t^n$ and showed that the
link is a compound torus knot formed by taking torus knot on a
small tube about the first
torus knot and iterating this procedure.
He also showed that only a finite number of terms (the characteristic
pairs) in the (possibly
infinite) power series parameterization of the
curve determined the topological type of the knot.
He continued by analyzing the case of curves with two branches.
Brauner concluded by computing the fundamental group of the complement of
these compound torus knots in terms of Wirtinger's generators and relations.

The next work in this area was done by Erich K\"{a}hler \cite{kahler}
in Leipzig, who remarks at the beginning of his paper that ``Obwohl die
betreffenden Fragen zum gr\"ossten Teil bereits von Herrn Brauner
beantwortet sind, habe ich mir erlaubt den Gegenstand auf dem etwas
anschaulicheren Wege...darzustellen.'' [Although this question has
been for the most part already answered by Mr. Branuner, I have allowed
myself to explain it in a somewhat clearer fashion.]

K\"ahler replaced 
Brauner's 
sphere, the boundary of the ``round'' four-ball $ \{
|x|^2 + |y|^2 \leq r^2\} $ in $\complex^2$ by the boundary of the
``rectangular'' four-ball  
$ \{ |x| \leq c_1\} \cap \{ |y| \leq c_2 \}$.
This is a simplification since a curve tangent to the $x$-axis (say)
intersects this boundary only in $\{|x| \leq c_1\} \cap \{|y| = c_2\}$,
one of its two sides ($c_1 << c_2$).
He noted that the two pieces of the boundary could be mapped easily into
three-space where they formed a decomposition into two solid tori.
He then looked at the curve
$y = ax^{m/n}$ and observed that the image of the intersection
of this curve with the boundary of the rectangular four-ball is obviously
a torus knot or link.
He then continued to obtain Brauner's results in easier fashion.

Thus the topological nature of the link could be computed from analytic data.
The converse result, that the characteristic pairs could be determined
from the
topology of the knot, was proved  
simultaneously by Zariski at Johns
Hopkins University and Werner Burau in K\"onigsberg.

Zariski \cite{zariski} started with a singular point of the curve $X$ and
again derived a presentation
of the local fundamental group of its complement.  He then found a polynomial
invariant $F(t)$ of this group which he later identified as the Alexander
polynomial of the knot, and showed that the first Betti number of the
$k$-fold branched cyclic cover of a punctured neighborhood of the origin in
$\complex^2$ with branch locus $X$ is the number of roots of
$F(t)$ which are $k$-th roots of unity.  (This was later recognized
to be a purely knot-theoretic result.)

%Also Zariski, commenting on the work of Brauner, mentions that:
%``The characteristic numbers 
%[the finite number of terms needed for the knot type] 
%of a singular branch determine completely the characteristic
%of the singularity in the sense of N\"other, i.e., as an aggregate of
%successive, infinitely near multiple points'', 
%thus making explicit the
%connection between topological and analytic data.

Burau \cite{burau-32}, on the other hand, 
used Alexander's
recent work to compute the Alexander
polynomial of compound torus knots.  He derived a
recursive formula for these polynomials and showed that they were all
distinct.  He later treated the case
when the polynomial had two branches at the origin, i.e. when the link 
had two components \cite{burau-34}.  

A survey of the above work was given later by Reeve \cite{reeve},
who also showed that the intersection number of two branches of a curve
at the origin equals the linking number in the three-sphere of their
corresponding knots.
He gave two proofs.
The first, following Lefschetz, notes that the algebraic
intersection multiplicity of the curves is their topological
intersection multiplicity, which is the linking number of their
boundaries.
The second proof uses Reidemeister's definition of linking number in
terms of the knot projection.

Now let us move forward in time to the present.
The computation of knot invariants of the link of a curve singularity becomes
increasingly messy as the number of branches of the curve increases.  A
diagrammatic method for these computations (for the Alexander
polynomial,the real Seifert form, the Jordan normal form of the
monodromy and so forth) has been developed in \cite{eisenbud-neumann}.

The link of a singularity of a curve has a global analogue, the {\em link
at infinity} $K_{\infinity}$ of a curve $X \subset \complex^2$, which is defined to be the intersection of
$X$ with a
sphere $\sphere^3_r$ of suitably large radius $r$.
Neumann has shown that if the curve is a regular fiber of its defining
equation (i. e. if the map is a locally trivial fibration near this
value), then the topological type of the curve is determined by the knot
type of $K_{\infinity} \subset \sphere^3_r$.
Also, Neumann and Rudolph have used these techniques to give topological
proofs of a result of Abhyankar and Moh (that up to algebraic automorphism,
the only embedding of $\complex$ in $\complex^2$ is the standard one)
and similar results of Zaidenberg and Lin
\cite{rudolph-82, neumann-rudolph, neumann-89}.

The knot type of the link of a singularity in higher dimensions has
received some attention; see for instance \cite{durfee-75}.

%%%%%%%%%%%%%%%%%%%%%%%%%%%%%%%%%%%

\section{Three-manifolds and singularities of surfaces}
\label{surface-section}

It is useful at this point to introduce some terminology.
An (affine) {\em algebraic variety} $X \subset \complex^m$ is the zero locus of a
collection of complex polynomials in $m$ variables.
If $X$ is a hypersurface, and hence the zero locus of a single polynomial $f(x_1, x_2, \dots,
x_m)$, then a point $p$ is {\em singular} if
$\del f /
\del x_1 = \dots = \del f / \del x_m = 0$ at $p$.
The set of nonsingular
points is a complex manifold of dimension $m-1$.
A point which is not singular is called {\em smooth}.
(The definition of singular point for arbitrary varieties can be found,
for example in \cite[Section 2]{milnor-spch}, and similar results hold.)
%The set of singular points of $X$ is proper subset
%outside of which $X$ is a smooth complex manifold. 

If $p \in X \subset \complex^m$, the {\em link} of $p$ in $X$ is
defined to be
$$K = X \cap \sphere^{2m-1}_{\epsilon}$$ 
where $\sphere^{2m-1}_{\epsilon}$ 
is a sphere of sufficiently small radius $\epsilon$ about $p$ in
$\complex^{m}$.  If $p$ is an isolated singularity of $X$,
then the link is a compact smooth real manifold of dimension one less
than the real dimension of $X$ at $p$.
Understanding the topology of the variety $X$ near
$p$ is the same as understanding the topology of $K$ and its embedding in the 
sphere; in fact, $X$ is locally homeomorphic to a cone on $K$ with vertex
$p$ \cite[2.20]{milnor-spch}.
(This fact is implicit in the work of Burau and K\"ahler, but not
explicitly stated.)
The {\em local fundamental group} of the singularity is the fundamental
group of the link.  
This is particularly interesting for an isolated singular point of an
algebraic surface (complex
dimension two) where the link is a three-manifold.

Some time elapsed before the topological investigation of curve
singularities chronicled in the first section was extended to higher
dimensions. 
In the early 1960's the following result of Mumford 
confirmed a conjecture of Abhyankar \cite{mumford} (see also the Bourbaki talk
of Hirzebruch \cite{hirzebruch-63}):

\begin{theorem}
If $p$ is a normal point of a complex surface $X$ with trivial local
fundamental group, then $p$ is a smooth point of $X$.
\end{theorem}

The condition ``normal'' comes from the algebraic side of algebraic
geometry; in particular it implies that the singularity is isolated and
that its link is a connected space.

%In particular, if the link is the three-sphere, then its embedding in
%$\sphere^{2N-1}_\epsilon$ is the standard embedding. PROOF

He proved this theorem by resolving the singularity, a technique which
in the case of surfaces is old and essentially algorithmic.
The process of resolution removes the singular point $p$ from $X$ and replaces
it by a collection of smooth transversally-intersecting complex curves
$E_1, \dots, E_r$ so that the new space $\tilde{X}$
is smooth.  

He showed that
the link could be obtained from the curves $E_i$  by a process called
{\em plumbing}:  
The tubular neighborhood of $E_i$ in $\tilde{X}$ is identified with
a 2-disk bundle over the curve $E_i$.  If $E_i$ and $E_j$ intersect in
a point $q \in \tilde{X}$, the two-disk bundles over $E_i$ and $E_j$ are glued
together by identifying a fiber over $q$ in one with a disk in the
base centered at $q$ in the other.  This makes a manifold with
corners. If the corners are smoothed (so that the result looks rather
like an plumbing elbow joint), the boundary is
diffeomorphic to the link.

The {\em graph} of a resolution of a normal singularity of
an algebraic surface is as follows:  The i-th
vertex corresponds to the curve $E_i$, labelled by the genus of $E_i$ and
the self-intersection $E_i \cdot E_i$.  The i-th and j-th vertices are
joined by an edge if $E_i \cdot E_j \neq 0$, and the edges are weighted
by $E_i \cdot E_j$.  The resolution graph thus determines the topological
type of the link.

Mumford used Van Kampen's theorem and the plumbing description of the link
to give a presentation of the local
fundamental group of the singularity and thus prove the theorem.

The local fundamental group of a singularity of an algebraic
surface turned out to be a a useful way to classify these singularities.
For instance, Brieskorn \cite{brieskorn-68}, using earlier work of Prill,
showed that if the local fundamental group is finite, then the variety
$X$ is locally isomorphic to a quotient $\complex^2 / G$, where $G$ is
one of the well-known
finite subgroups of $GL(2,\complex)$. 
He listed all such subgroups $G$ together with the resolution graph of the
minimal resolution of the corresponding singularity $\complex^2 / G$.

Wagreich \cite{wagreich-72}, inspired by work of Orlik
\cite{orlik-70},  used Mumford's presentation to find all
singularities with nilpotent or solvable local fundamental group.
Thus the local fundamental group became closely connected with the local
analytic structure of the singularity.

Neumann showed that the topology of the link $K$ determines the graph of
the minimal resolution of the singularity.
In fact, he showed
that $\pi_1(K)$ determines this graph, except in a small number
of cases \cite{neumann-81} .

Mumford's techniques in a global setting appeared later in work of
Ramanujam \cite{ramanujam}:

\begin{theorem}
A smooth complex algebraic surface which is contractible and simply
connected at infinity is algebraically isomorphic to $\complex^2$.
\end{theorem}

Ramanujam showed this by compactifying the surface by a divisor with
normal crossings, and then using the topological conditions to show that
this divisor could be contracted to a projective line.
He also showed that the condition of simple connectivity at infinity was
essential by producing an example of a smooth affine rational surface
$X$ which is contractible but not algebraically isomorphic to
$\complex^2$.  In fact, the intersection of $X$ with a sufficiently
large sphere is a homology three-sphere but not a homotopy
three-sphere.

Ramanujam's result implies that the only complex algebraic
structure on $\real^4$ is the standard one on $\complex^2$, so that
there are no ``exotic'' algebraic structures on the complex plane.
The search for exotic algebraic structures thus continued in higher
dimensions.
Ramanujam remarked that the three-fold $X
\times \complex$ is diffeomorphic to
$\complex^3$ by the h-cobordism theorem. 
A cancellation theorem proved later had the corollary that $X
\times \complex$ is not algebraically
isomorphic to $\complex^3$.
Hence there is an exotic algebraic structure on $\complex^3$.
Much work followed in this area; for the current state of affairs one
can consult \cite{zaidenberg-93}, for example.

%%%%%%%%%%%%%%%%%%%%%%%%%%%%%%%%%%%%%%

\section{Exotic spheres}
\label{exotic-spheres-section}

Brieskorn, who was spending the academic year 1965-66 at the
Massachusetts Institute of Technology, investigated whether
Mumford's theorem extended to higher dimensions.
On September 28, 1965, he wrote in a letter to his doctoral advisor Hirzebruch
that he had examined the three-dimensional variety 
$$x_1^2 + x_2^2 + x_3^2 + x_4^3 = 0$$
and its singularity at the origin.
He
explicitly calculated a resolution of the singular point, then used van
Kampen's theorem to show that the link $K$ of this singularity is
simply-connected and the Mayer-Vietoris sequence to show that $K$ is a 
homology $5$-sphere.  He concluded, using Smale's recent solution of
the Poincar\'e conjecture in higher dimensions, that $K$ is
homeomorphic to $\sphere^5$.
Hence Mumford's result did not extend to higher dimensions.

According to Hirzebruch \cite[C38]{hirzebruch-works}, ``Dieser Brief von
Brieskorn war eine grosse \"Uberraschung'' [This letter from Brieskorn was a great surprise].
Later letters followed with more squared terms added
to the equation above.
Brieskorn's final result appeared in
\cite{brieskorn-66-1}: 
For odd $n \geq 3$, the link at the origin of 
\begin{equation}
\label{kervaire-link}
x_0^3 + x_1^2 + x_2^2 + \dots + x_n^2 = 0 
\end{equation}
is homeomorphic to the sphere $S^{2n-1}$.  

The attention then shifted to the differentiable structure on this link.
To describe the next events, we first need to recall the situation with
non-standard or ``exotic'' differentiable structures on spheres.
The first exotic sphere, a differentiable structure on $\sphere^7$ which
is not diffeomorphic to the standard structure, had been discovered only ten years earlier by
Milnor.  Further investigations followed by Kervaire and Milnor
\cite{kervaire-milnor}.  By Smale's solution to the
higher-dimensional Poincar\'e conjecture, it was sufficient to
look at the set $\Theta_m$ of homotopy $m$-spheres (manifolds homotopy
equivalent to the standard sphere $S^m$).  The set $\Theta_m$ is an
abelian group under connected sum, and Kervaire and Milnor showed that this
group is finite ($m \neq 3$).  

They also looked at the subgroup $bP_{m+1} \subset
\Theta_m$ of homotopy spheres which are boundaries of parallelizable
manifolds (manifolds with trivial tangent bundle), and showed that
$bP_{m+1}$ is trivial for $m$ even, and finite cyclic for $m \neq 3$ odd.  

Its order could be computed as follows: If $n$ is odd, the group
$bP_{2n}$ has order one or two.  It is generated by the Kervaire
sphere which is the boundary of the manifold constructed by plumbing
two copies of the tangent disk bundle to $S^n$.  The Kervaire sphere may or
may not be diffeomorphic to the standard sphere; the first non-trivial
group is $bP_{10}$.  If $\Sigma \in bP_{2n}$ is the boundary of an
$(n-1)$-connected parallelizable $2n$-manifold $M$, whether $\Sigma$ is
diffeomorphic to the standard sphere or the Kervaire sphere depends on
the Arf invariant of a geometrically-defined quadratic form on $M$.

If $n \geq 4$ is even, the order of $bP_{2n}$ can be calculated in terms
of Bernoulli numbers. For example, there are 28 homotopy seven-spheres
in $\Theta_7 = bP_{8}$.  Also, the order of $\Sigma \in bP_{2n}$ can be
calculated in terms of the signature of the intersection pairing on
$H^{n}(M)$.

\begin{figure}
\postscript{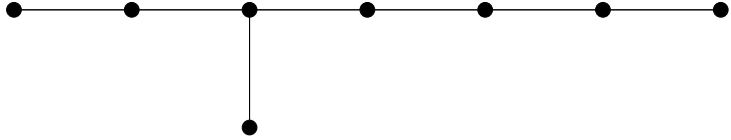}{0.7}
\caption{The $E_8$ graph}
\label{e8-graph}
\end{figure}

The construction of a generator of $bP_{2n}$ for $n \geq 4$ even is
once again bound up with singularity theory.  In a preprint
\cite{milnor-59} of January, 1959, Milnor had constructed a generator by
plumbing according to an even unimodular matrix of rank and index eight.
This matrix was not the well-known one associated to the $E_8$ graph
(Figure \ref{e8-graph}), 
though, since its graph had a cycle.  He then added a two-handle to make the
boundary simply-connected and hence a homotopy sphere.  Hirzebruch,
however, was familiar with the $E_8$ matrix from his work on resolution
of singularities of surfaces.  He constructed a generator of the group
$bP_{2n}$ 
by plumbing copies of the tangent disk bundle to $S^n$ according to
the $E_8$ graph.  (For more details, see
\cite[C30]{hirzebruch-works}, \cite{hirzebruch-63, hirzebruch-66, milnor-64}.)

At the same time in the fall of 1965 that Hirzebruch was receiving the
letters from Brieskorn, he also received a letter from Klaus J\"anich,
another of his doctoral students, who was spending the year 1965-66 at
Cornell.  J\"anich described his work on $(2n-1)$-dimensional
$O(n)$-manifolds (manifolds with an action of the orthogonal group).  In
fact, he had classified $O(n)$-manifolds whose action had just two orbit
types with isotropy groups $O(n-1)$ and $O(n-2)$, in particular showing
that they were in one-to-one correspondence with the non-negative
integers.  (These results were also obtained by W. C. Hsiang and
W. Y. Hsiang.)

\begin{figure}
\postscript{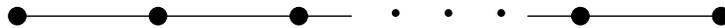}{0.7}
\caption{The $A_{k}$ graph ($k$ vertices)}
\label{a-graph}
\end{figure}

Hirzebruch noticed the connection between the research efforts of his two
students and showed that the link of
\begin{equation}
\label{d-variety}
x_0^d + x_1^2 + \dots + x_n^2 = 0
\end{equation}
for $d \geq 2$, $n \geq 2$ is an $O(n)$-manifold as above with invariant $d$, the action being
given by the obvious one on the last $n$ coordinates.
Since the boundary of the manifold constructed by plumbing copies of the
tangent disk bundle of the $n$-sphere according to the $A_{d-1}$ tree
(Figure \ref{a-graph})
also is an $O(n)$-manifold as above with invariant $d$,
these manifolds are identical.
Thus the link
of the singularity (\ref{kervaire-link}) is the
$(2n-1)$-dimensional Kervaire sphere;
in particular for $n = 5$ it is an exotic $9$-sphere.

These results were described in a manuscript
``$O(n)$-Mannigfaltigkeiten, exotische Sph\"aren, kuriose Involutionen''
of March 1966.
(This was not published, since it was supplanted by Hirzebruch's Bourbaki talk
\cite{hirzebruch-66}, and the detailed lecture notes
\cite{hirzebruch-mayer} from his course in the winter semester 1966/67
at the University of Bonn.)
In a letter \cite[C39]{hirzebruch-works} of March 29, 1966, Brieskorn reacted to the manuscript with
``Klaus J\"anich und ich hatten von diesem Zusammenhang unserer Arbeiten
nichts bemerkt, und ich war vor Freude ganz ausser mir, wie Sie nun die
Dinge zusammengebracht haben.  Ein sch\"oneres Zusammenspiel von Lehrern
und Sch\"ulern--wenn ich das so sagen darf--kann man sich doch wirklich
nicht denken.''
[Klaus J\"anich and I had not noticed this connection between our work,
and I was beside myself with joy to see how you had brought these
together.  A more beautiful cooperation of student and pupil can one
hardly imagine, if I may say so myself.]

At this time the varieties
\begin{equation}
\label{brieskorn-variety}
x_0^{a_0} + x_1^{a_1} + \dots + x_n^{a_n} = 0
\end{equation}
($a_i \geq 2$) started to receive attention; they are now
called ``Brieskorn
varieties'', probably due to the influence of a chapter heading in
Milnor's book \cite{milnor-spch}, although they were first examined in
this context by
Pham and Milnor as well.
The corresponding $(2n-1)$-dimensional links 
$$K(a_0, a_1, \dots, a_n) = \{x_0^{a_0} + x_1^{a_1} + \dots + x_n^{a_n}
= 0 \} \cap \sphere^{2n+1} $$
where $\sphere^{2n+1}$ is a sphere about the origin, are usually
called ``Brieskorn manifolds''.
(The radius of the sphere can be arbitrary since the equation is
weighted homogeneous.)

Milnor, who was in Princeton, sent a letter in April of 1966 to John
Nash at MIT describing a simple conjecture as to when $K(a_0, a_1,
\dots, a_n)$ is a homotopy sphere: Let $\Gamma(a_0, a_1,
\dots, a_n)$ be the graph with $n+1$ vertices labeled $0, 1, \dots, ,n$ and with two
vertices $i$ and $j$ joined by an edge if the greatest common divisor
$(a_i, a_j)$ is bigger than 1.

\begin{figure}
\postscript{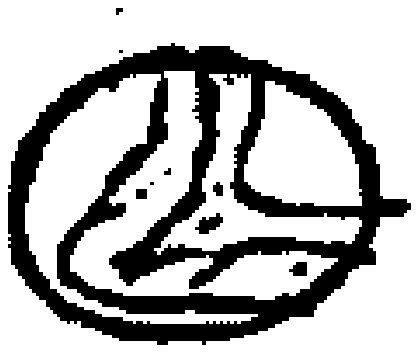}{0.2}
\caption{Milnor's sketch}
\label{figure4}
\end{figure}

\begin{conjecture}
For $n \geq 3$, 
the link $K(a_0, a_1, \dots, a_n)$ is a homotopy $(2n-1)$-sphere if and only if
the graph $\Gamma(a_0, a_1, \dots, a_n)$ has 
\begin{itemize}
\item at least two isolated
points, or 
\item one isolated point and at least one connected component $\Gamma'$
with an odd number of vertices such that the g.c.d. $(a_i, a_j) = 2$ for all
$i \neq j \in \Gamma'$.
\end{itemize}
\end{conjecture} 
In the corner of Milnor's letter was a sketch (Figure \ref{figure4}) which has now
become the standard picture of a map with an isolated critical point.

Brieskorn then chanced upon an article of Pham \cite{pham} which dealt
with exactly the variety (\ref{brieskorn-variety}) above.
In fact,
Pham was interested in calculating the 
ramification of certain integrals encountered in the interaction of
elementary particles in theoretical physics.
To do this he needed to
generalize the
Picard-Lefschetz formulas, so let us recall these.

Picard-Lefschetz theory can be summarized as follows (see for example
\cite[2.1]{encyc}):
Let
$$X_t = \{x_0^{2} + x_1^{2} + \dots + x_n^{2} = t \} \subset \complex^{n+1}$$
%so that $X_0$ is a cone, and $X_t$ for $t \neq 0$ is smooth.
($n \geq 1$).  Then
\begin{enumerate}
\item The smooth variety $X_t$ for $t \neq 0$ is homotopy equivalent to an
$n$-sphere $S^n$.  (In fact, it is diffeomorphic to the tangent bundle
to $S^n$.) 
\item The homology class of this $n$-sphere generates the
kernel of the degeneration map $H_n(X_t) \to H_n(X_0)$, hence its name of {\em vanishing
cycle}.
\item The self-intersection of the vanishing cycle is 2 if $n \congruent
0$ mod 4, $-2$ if $n \congruent 2$ (mod 4) and 0 if $n \congruent
1, 3$ (mod 4).
\item Starting at $t=1$ in the complex plane, traveling once
counterclockwise about 
the origin and returning to the starting point induces a smooth map
called the {\em monodromy} of
$X_1$ to itself.  It is well-defined up to isotopy.
Picard-Lefschetz theory gives a description of this map.
For example, if $n=1$ it is
a Dehn twist about the one-dimensional vanishing cycle.
Picard-Lefschetz theory also describes the induced maps
$H_n(X_1) \to H_n(X_1)$ and $H_n(X_1, \boundary X_1) \to H_n(X_1)$. 
\end{enumerate}

Pham generalized this situation to the case
$$X_t = \{x_0^{a_0} + x_1^{a_1} + \dots + x_n^{a_n}=t \} \subset
\complex^{n+1} $$
and found
\begin{enumerate}
\item The smooth variety $X_t$ for $t \neq 0$ is homotopy equivalent to a
bouquet $S^n \vee S^n \vee \dots \vee S^n$ of $(a_0 -1)(a_1 -1) \dots (a_n
-1)$ $n$-spheres.
(This was shown by retracting $X_t$ to a join $Z_{a_0} * Z_{a_1} * \dots
* Z_{a_n}$ where $Z_{k}$ denotes $k$ disjoint points.)
\item The homology classes of these $n$-spheres generate the kernel of
the map $H_n(X_t) \to H_n(X_0)$.
\item An explicit calculation of the intersection pairing on $H_n(X_t)$.
\item An explicit calculation of the monodromy action on $H_n(X_t)$. (This is induced
by rotating each set of points $Z_k$.)
\end{enumerate}

The article of Pham provided exactly the information Brieskorn needed.
(He remarks \cite{brieskorn-66-2} that ``F\"ur den Beweis von [diesen] Aussagen sind
jedoch gewisse Rechnungen erforderlich, f\"ur die gegenw\"artig keine
allgemein brauchbare Methode verf\"ugbar ist.  F\"ur den Fall der
$K(a_0, a_1, \dots, a_n)$ sind diese Rechnungen aber s\"amtlich in einem
vor kurzem erschienen Artikel von Pham enthalten, und nur die Arbeit von
Pham erm\"oglicht den so m\"uhelosen Beweis unserer Resultate.''
[Certain calculations, for which there are no general methods at this time,
are necessary for the proof of these results.  In the case of
$K(a_0, a_1, \dots, a_n)$, however, these calculations are contained in
an article of Pham which just appeared, and it is only Pham's work which
makes possible such an effortless proof of our results.])
Brieskorn used it to prove a 
conjecture of Milnor from the preprint \cite{milnor-66} about the
characteristic polynomial of the monodromy \cite[Lemma
4]{brieskorn-66-2}, \cite[Theorem 9.1]{milnor-spch}.
He then used this to prove 
the conjecture above \cite[Satz 1]{brieskorn-66-2},
\cite[14.5]{hirzebruch-mayer},\cite[Section 2]{hirzebruch-66}.

Brieskorn also noted that the link $K(a_0, a_1, \dots, a_n)$,
which is defined as $X_0 \cap \sphere^{2n+1}$, is
diffeomorphic to $X_t \cap \sphere^{2n+1}$ for small $t \neq
0$. 
This is the boundary of the smooth $(n-1)$-connected manifold $X_t \cap \ball^{2n+2}$;
it is
parallelizable since it has trivial normal bundle.
Hence $K(a_0, a_1, \dots, a_n) \in bP_{2n}$.
The information in Pham's paper about the intersection form also led to a
formula (derived by Hirzebruch) for the signature of $X_t \cap \ball^{2n+2}$.
Brieskorn concluded that the link of
\begin{equation}
\label{exotic-7}
x_0^{6k-1} + x_1^3 + x_2^2 + x_3^2 + \dots + x_n^2 = 0
\end{equation}
for even $n \geq 4$ is $k$ times the Milnor generator of $bP_{2n-1}$
\cite{brieskorn-66-2, hirzebruch-mayer}.

Through a preprint of Milnor \cite{milnor-66}, Brieskorn also learned of
a recent result of Levine \cite{levine}
which showed how to compute the Arf invariant needed to recognize
whether a link is the Kervaire sphere in terms of the higher-dimensional
Alexander polynomial of the knot.  The Alexander polynomial for
fibered knots is the same as the characteristic polynomial of the
monodromy on $H_n(F)$.  Hence Brieskorn was
able to
show \cite[Satz 2]{brieskorn-66-2} that the link of
$$x_0^d + x_1^2 + x_2^2 + \dots + x_n^2 = 0 $$
for $n \geq 3$ odd
is the
standard sphere if $d \congruent \pm 1 $ mod 8, and the Kervaire sphere
if $d \congruent \pm 3 $ mod 8, thus providing 
another proof of Hirzebruch's result that the link of the
singularity (\ref{kervaire-link}) is the Kervaire sphere.

The explicit representation of all the elements of $bP_{2n}$ by links of
simple algebraic equations was rather surprising.  
It provided another way of thinking about these
exotic spheres and led to various topological applications.

For example, Kuiper \cite{kuiper} used them to obtain algebraic
equations for all non-smoothable piecewise-linear manifolds of dimension
eight.
(PL manifolds of dimension less than eight are
smoothable.) 
In fact, he started with the complex four-dimensional variety given by
Equation (\ref{exotic-7}) above with $n = 4$.
This has a single isolated singularity at the origin.
Its completion in projective space has singularities on the hyperplane
at infinity, but adding terms of higher order to the equation eliminates
these while keeping (analytically) the same singularity at the origin.
This can be triangulated, giving a combinatorial eight-manifold which is
smoothable except possibly at the origin.
Since obstructions to smoothing are in one-to-one correspondence with
the 28 elements of $bP_8$, the construction is finished.

Also, the high symmetry of the variety given by Equation
(\ref{d-variety}) allowed the construction of many interesting group
actions on spheres, both standard and exotic \cite[Section
4]{hirzebruch-66}, \cite[Section 15]{hirzebruch-mayer}.  The actions are
the obvious ones: The cyclic group of order $d$ acts by roots of unity
on the first coordinate, and there is an involution acting on (any
subset of) the remaining coordinates by taking a variable to its
negative.

%%%%%%%%%%%%%%%%%%%%%%%%%%%%%%%%%%%%%%%

\section{The Milnor fibration}
\label{section-milnor-fibration}

About the same time as the above events were happening, Milnor proved a
fibration theorem which turned out to be fundamental for much subsequent work.  
This theorem
together with its consequences first appeared in the unpublished
preprint \cite{milnor-66}, which
dealt exclusively with isolated
singularities.  
(A full account of this work was later published in the book
\cite{milnor-spch}, where the results were generalized to 
non-isolated singularities.  The earlier and somewhat simpler ideas can be found at the end of
Section 5 of the book.) 

Let $f(x_0, x_1, \dots, x_n)$ for $n \geq 2$ be a complex polynomial with $f(0, \dots,
0) = 0$ and an isolated critical point at the origin.
Let $\sphere^{2n+1}_{\epsilon}$ be a sphere of suitably small radius
$\epsilon$ about the
origin in $\complex^{n+1}$.
As before, let $K = \{f(x_0, x_1, \dots, x_n) = 0 \} \cap
\sphere^{2n+1}_{\epsilon}$ be the link of $f = 0$ at the origin.  
The
main result of the preprint is the following {\em fibration theorem}:

\begin{theorem} 
The complement of an open tubular neighborhood of the link $K$ in
$\sphere^{2n+1}_{\epsilon}$ is the total
space of a smooth fiber bundle over the circle $\sphere^1$.
The fiber $F$ has boundary diffeomorphic to $K$.
\end{theorem}

The idea of the proof is as follows:
If $\ball^{2n+2}_{\epsilon}$ is the ball of radius $\epsilon$
about the origin and $\delta > 0$ is suitably small, then
$$ f: f^{-1}(\sphere^1_\delta) \cap \ball_{\epsilon}^{2n+2} 
\to \sphere^1_\delta  $$ 
is a clearly a
smooth fiber bundle with fiber
$$ F' = \{ f(x_0, x_1, \dots, x_n) =\delta \} \cap \ball^{2n+2}_{\epsilon}$$ 
The total space of this fibration is then pushed out to the sphere
$\sphere^{2n+1}_{\epsilon}$ along the trajectories $p(t)$ of a
suitably-constructed vector field.  This vector
field is constructed with the property that $|p(t)|$ is increasing along
a trajectory, so that points eventually reach the sphere, and with the
property that the argument of $f(p(t))$ is constant and $|f(p(t))|$ is
increasing, so that the images of points in $\complex$ travel out on
rays from the origin.
Thus Milnor's proof shows that $F$ is diffeomorphic
to $F'$.
The proof also shows that $F$ is parallelizable, since $F'$ has trivial normal
bundle. 

The fiber $F$ is now called 
 the {\em Milnor fiber}.
He then gives some facts which lead to the topological type of
the fiber $F$ and the link $K$:

\medskip

\noindent (a) \  The pair $(F, \boundary F)$ is
$(n-1)$-connected.

\medskip

\noindent (b) \  The fiber $F$ has the homotopy type of a cell
complex of dimension $\leq n$.  In fact, it is built from the $2n$-disk
by attaching handles of index $\leq n$.

\medskip

These assertions follow from Morse theory.  In fact, in a lecture at
Princeton in 1957 (which was never published), Thom described an
approach to the Lefschetz hyperplane theorems which was based on Morse
theory.  Thom's approach then inspired Andreotti and
Frankel \cite{andreotti-frankel-59} (see also \cite[Section
7]{milnor-mt}) to give another proof of Lefschetz's first hyperplane
theorem which used Morse theory, but in a different way:  The key
observation is that given a $n$-dimensional complex variety $X \subset \complex^m$ and a (suitably general) point $p \in \complex^m -
X$, then the function on $X$ defined by $|x - p |^2$ for $x \in X$ 
has non-degenerate
critical points of Morse index $\leq n$.
Thus $H_k(X) = 0$ for $k > n$, which is equivalent to Lefschetz's first
hyperplane theorem.

Assertion (b) follows since the function $|x|$ (or a slight
perturbation of it) restricted to $F'$ has critical points of index
$\leq n$, and Assertion (a) follows
since the function $-|x|^2$ on $F'$ has critical points of index
$\geq n$.

By (b), the complement $\sphere^{2n+1}_\epsilon - F$ has the same homotopy
groups as $\sphere^{2n+1}_\epsilon$ through dimension $n-1$.  Thus:

\medskip

\noindent (c)  The complement $\sphere^{2n+1}_\epsilon - F$ is
$(n-1)$-connected.

\medskip

By the Fibration Theorem, $\sphere^{2n+1}_\epsilon - F$ is homotopy equivalent
to $F$.  Thus

\begin{proposition}  The fiber $F$ has the homotopy type of a bouquet $S^n
\vee \dots \vee S^n$ of spheres.
\end{proposition}

Fact (a) and the above proposition combined with the long exact sequence
of a pair show the following:

\begin{proposition}
The link $K$ is $(n-2)$-connected.
\end{proposition}

%The existence of a particular fibration of the sphere given by the
%fibration theorem was used by Lawson \cite{lawson} and Durfee
%\cite{durfee-72} to construct
%codimension-one foliations of odd-dimensional spheres. 

Milnor used the notation $\mu$ for 
the number of spheres in the bouquet of the first proposition and
called it the ``multiplicity'' since it is the multiplicity of the
gradient map of $f$.
However, $\mu$ quickly became known as
the {\em Milnor number}. 
The Milnor number has played a central role in the study of
singularities.
One reason is that it has analytical as well as a topological
descriptions, for example: 
$$\mu = dim_{\complex} \complex\{x_0, x_1, \dots, x_n\}/ (\del f/\del x_0, \del
f/\del x_1, \dots, \del f/\del x_n ) $$
the (vector space) dimension of the ring of power series in $n+1$ variables divided by
the Jacobian ideal of the function
(see, for example, \cite{orlik-76}).

The fact that the Milnor number can be expressed in different ways is
extremely useful.
For example, the topological interpretation of $\mu$ was used by Le and
Ramanujam to prove a result which became basic to the study of
equisingularity:  Suppose that $n \neq 2$.  
If a family of functions $f_t: \complex^{n+1} \to \complex$ depending on $t$
with isolated
critical points has constant
Milnor number, then the differentiable type of the
Milnor fibration of $f_t$ is independent of $t$.
The proof uses the topological interpretation of $\mu$ to produce a h-cobordism
which is thus a product cobordism; hence the restriction $n \neq 2$
\cite{le-ramanujam}. 

Results similar to the Fibration Theorem and the two
propositions have now been obtained in many different situations:
complete intersections, functions on arbitrary varieties, polynomials
with non-isolated critical points, critical points of polynomials at
infinity, and so forth; references to these results can be found in the
books and conference proceedings cited at the beginning of this article.
Also, there are now many different techniques for computing the Milnor
number $\mu = \mbox{ \ rank \ } H_n(F)$,
the characteristic polynomial of the monodromy $H_n(F) \to H_n(F)$, and
the intersection pairing on $H_n(F)$.

The characteristic polynomial of the monodromy turned out to be
cyclotomic, and a variety of proofs have appeared of this
important fact: the geometric proof of Landman, geometric
proofs of Clemens and Deligne-Grothendieck based on resolving the
singularity, proofs based on the Picard-Fuchs equation by Breiskorn,
Deligne and Katz, and analytic proofs using the classifying space for
Hodge structures by Borel and Schmid.  For a summary of these and
the appropriate references, see \cite{griffiths-73}.

The Milnor number appears in another situation.  To describe this we
first return to Thom's original observation in his 1957 lecture, as recorded in
\cite{andreotti-frankel-69}: Given an $n$-dimensional
complex variety $X$ in affine space and a suitably general linear
function $f: X \to \complex$, then $|f|^2$ has non-degenerate critical points of Morse
index exactly $n$ (except for the absolute minimum).  This result is
easily proved by writing the function in local coordinates.  It forms
the basis of Andreotti and Frankel's proof of the second hyperplane
theorem of Lefschetz, which says that the kernel of the map on $H_{n-1}$
from a hyperplane section of an $n$-dimensional projective
variety to the variety is generated by vanishing cycles.

Thom's original observation was applied in the local context of singularities,
where it leads to a basic result in the subject of polar curves
relating the Milnor number of a singularity and a plane section.
This result has both topological and analytic formulations
\cite[p. 317]{teissier-73}, \cite{le-73}.

%%%%%%%%%%%%%%%%%%%%%%%%%%%%%%%

\section{Brieskorn three-manifolds}

The Brieskorn three-manifolds $K(a_0, a_1, a_2)$, the link of 
$$x_0^{a_0} +
x_1^{a_1} + x_2^{a_2} = 0$$
 at the origin, 
have
provided examples figuring in many topological
investigations. 
For example, 
the local fundamental group of these singularities has proved interesting.
As mentioned in Section \ref{surface-section}, the surface singularities whose link have finite
fundamental group are exactly the quotient singularities.
If the surface is embedded in codimension one, and is hence the zero
locus of a polynomial $f(x_0, x_1, x_2)$, then these singularities are
the well-known {\em simple singularities}:
$$A_k:  x_0^{k+1} + x_1^2 + x_2^2 = 0 \ \  (k \geq 1)$$
$$D_k:  x_0^{k-1} + x_0x_1^2 + x_2^2 = 0 \ \ (k \geq 4)$$ 
$$E_6:  x_0^4 + x_1^3 + x_2^2 = 0 $$
$$E_7:  x_0^3 + x_0 x_1^3 + x_2^2 = 0$$
$$E_8:  x_0^5 + x_1^3 + x_2^2 = 0$$
These equations have appeared, and continue to appear, in many seemingly
unrelated contexts \cite{durfee-79}.
For example, V. Arnold showed that they are the germs of functions whose
equivalence classes under change of coordinate in the domain have no
moduli \cite{agv-1}.

More general than Brieskorn polynomials is the class of weighted homogeneous polynomials:
A polynomial $f(x_0, x_1, \dots, x_n)$ is {\em weighted homogeneous} if
there are positive rational numbers $a_0, a_1, \dots, a_n$ such that 
$$f(c^{1/a_0}x_0, c^{1/a_1}x_1, \dots , c^{1/a_n}x_n) = c f(x_0, x_1,
\dots, x_n)$$
for all complex numbers $c$.
(Weighted homogeneous polynomials probably first made their appearance in
singularity theory in the book of Milnor \cite{milnor-spch}.) 
Brieskorn singularities are weighted homogeneous, with weights
exactly the exponents.

The simple singularities are weighted homogeneous.  Milnor
\cite[p. 80]{milnor-spch} noted that their weights $(a_0, a_1, a_2)$
satisfy the inequality $1/a_0 + 1/a_1 + 1/a_2 > 1$.  He also remarked
that the links of the {\em simple elliptic singularities}
$$\tilde{E}_6:  x_0^3 + x_1^3 + x_2^3=0$$
$$\tilde{E}_7:  x_0^2 + x_1^4 + x_2^4=0$$
$$\tilde{E}_6:  x_0^2 + x_1^3 + x_2^6=0$$
have infinite nilpotent fundamental group.
In this case, the sum of the reciprocals of the weights is 1.
He conjectured that if $1/a_0 + 1/a_1 + 1/a_2 \leq 1$, then the
corresponding link had infinite fundamental group, and that this group
was nilpotent exactly when  $1/a_0 + 1/a_1 + 1/a_2 = 1$.

This conjecture was proved by Orlik \cite{orlik-70}.
In fact, he and Wagreich
\cite{orlik-wagreich} had already found an explicit form of a resolution
for weighted homogeneous singularities using topological methods based
on the existence of a $\complex^*$ action, following earlier work by
Hirzebruch and J\"anich.
They also noted that these links were Seifert manifolds
\cite{seifert-32} and hence could use Seifert's work as well as
earlier work by Orlik and others.

Topologists were interested in the question of which homology
three-spheres bound contractible four-manifolds (c.f. \cite[Problem
4.2]{kirby-78}).
In fact, topological analogues (contractible four-manifolds which are
not simply-connected at infinity) of the example of Ramanujam in Section
\ref{surface-section} (a contractible complex surface which is not simply
connected at infinity) had been found some ten
years earlier by Mazur \cite{Mazur} and Poenaru \cite{Poenaru}.
As Mazur remarks, these examples provide a method of constructing many
examples of odd topological phenomena.

It was known (see Milnor's conjecture in Section \ref{exotic-spheres-section})
that $K(a_1,a_2,a_3)$ is a homology three-sphere exactly
when the integers $a_1,a_2,a_3$ are pairwise relatively prime.
(As Milnor remarks in \cite{milnor-75}, 
this result in this context of Seifert fiber spaces is already in
\cite{seifert-32}.) 
Links of Brieskorn singularities were particularly easy to study, since
a resolution of the singularity exhibited the link as the boundary of a
four-manifold, and data from the resolution provided a
plumbing description of this manifold which then could be manipulated to
eventually get a contractible manifold.
For example, Casson and Harer \cite{casson-harer} showed that the
Brieskorn manifolds $K(2,3,13)$, $K(2,5,7)$ and
$K(3,4,5)$ are boundaries of contractible four-manifolds.
Much has now happened in this area as can be seen in Kirby's update of
his problem list \cite{kirby-book}.

Brieskorn three-manifolds and their generalizations also provided
interesting examples of manifolds with a ``geometric structure''.
Klein proved long ago that the links of the simple singularities listed
above are of the
form $S^3/\Gamma$, the quotient of the group of unit quarternions by
a discrete subgroup.

Milnor \cite[Section 8]{milnor-75} proved by a round-about method that
the links of the simple elliptic singularities are quotients of the
Heisenberg group by a discrete subgroups.  He then showed that the links
of Brieskorn singularities with $1/a_0 + 1/a_1 + 1/a_2 \leq 1$ are
quotients of the universal cover of $SL(2,\real)$ by discrete
subgroups. (Similar results were obtained at the same time by
Dolgachev.)

Thus many links admitted a locally homogeneous (any two points have
isometric neighborhoods) Riemannian metric and
hence provided nice examples of Thurston's eight geometries \cite{thurston}.
These results were extended by Neumann \cite{pspm-40}.
Later he and Scherk \cite{neumann-scherk} found a more natural way of
describing the connection between the geometry on the link and the
complex analytic structure of the singularity in terms of locally
homogeneous non-degenerate CR structures.

The three-dimensional Brieskorn manifolds have also been central
examples in the study of the group $\Theta^H_3$ of homology
three-spheres.  This group is bound up with the question of whether
topological manifolds can be triangulated.  It was originally thought
that this group might just have two elements.  However, techniques from
gauge theory were used to show that it is actually infinite and even
infinitely generated.  In particular the elements $K(2,3,6k-1)$ for $k
\geq 1$ have infinite order in this group, and are linearly independent.
Brieskorn manifolds appear in this context because the
corresponding singularities have easily computable resolutions, and
hence the three-manifolds are boundaries of plumbed four-manifolds upon
which explicit surgeries can be performed \cite{fintushel-stern}.

Also, the Casson invariant of some types of links of surface singularities in
codimension one (including Brieskorn singularities) was proved to be
$1/8$ of the signature of the Milnor fiber  \cite{neumann-wahl}.

\section{Other developments}

This last section recounts two developments which occurred outside the
main stream of events as recounted in the previous sections.
They are both applications of topology to algebraic geometry.
The first is a
theorem of Sullivan \cite{sullivan}: 

\begin{theorem}
If $K$ is the link of a point in a complex algebraic variety, then the
Euler characteristic of $K$ is zero.
\end{theorem}

If the point is smooth or an isolated singular point, then the link is a
compact manifold of odd dimension and hence has Euler characteristic
zero.  The surprising feature of this result is that is should be true
for non-isolated singularities as well.

Sullivan discovered this result during his study of
combinatorial Stiefel-Whitney classes.
He recounts that initially it was clear to him that this result
was true in dimensions one and two.  He then asked Deligne if he knew of
any counterexamples in higher dimensions, but the latter replied
``almost immediately'' with a proof based on resolving the singularity.
Sullivan then deduced
this result in another fashion: Since complex varieties have a
stratification with only even-dimensional strata, the link has a
stratification with only odd-dimensional strata.  He then proved, by
induction on the number strata, that a compact stratified space with
only odd-dimensional strata has zero Euler characteristic.

Since real varieties are the fixed point set of the conjugation map
acting on their complexification, this result has the following
consequence for real varieties:

\begin{corollary}
If $K$ is the link of a point in a real algebraic variety, then the
Euler characteristic of $K$ is even.
\end{corollary}

Sullivan remarks that the result for complex varieties follows from
essentially ``dimensional considerations'', but that the corollary for
real varieties is however ``geometrically surprising''.
This result continues to form a basis for the investigation of the
topology of real varieties.

The second result is one of Thom.
Given a singularity of an arbitrary variety
$X_0 \subset \complex^m$, one can ask if it can 
be ``smoothed'' in its ambient space $\complex^m$ in the sense that it
can be made a fiber of a flat family $X_t \subset \complex^m$
whose fibers $X_t$ for small $t \neq 0$ are smooth. 
For example, a hypersurface singularity is smoothable in its
ambient space since it is the zero locus of a
polynomial and hence smoothed by nearby fibers of the polynomial.

The first example of a non-smoothable singularity was constructed by Thom (see \cite{hartshorne}).
In fact, Thom showed that the variety $X \subset \complex^6$ defined by
the cone on the Segre embedding $\projective^1 \times \projective^2 \to
\projective^5$ is not smoothable:
If it were, the link $K^7 \subset \sphere^{11}$ of $X \subset
\complex^6$ at the origin would be null cobordant (as a manifold with
complex normal bundle) in $\sphere^{11}$,
but it is not.  
This is proved by a computation with characteristic classes.
(The manifold $K^7$
is odd-dimensional and hence null
cobordant, but not in $\sphere^{11}$.)

%%%%%%%%%%%%%%%%%%%%%%%%%%%%%%%%%%%%%%%%%%%%%%%%%%

\bigskip

Acknowledgments: I thank F. Hirzebruch allowing me to use material from
a lecture in July 1996 at the Oberwolfach conference in honor of
Brieskorn's 60th birthday, and I thank both him and D. O'Shea for
comments on a preliminary version of this article.

\end{document}